\documentclass[12pt]{article}
 \usepackage{epsfig}
\usepackage{graphics}
 \newcommand{\hx}{\hat x}
 \newcommand{\tx}{\tilde x}
\newcommand{\Mo}{\mathcal{M}}

\usepackage{color}
\usepackage{epsfig,cmmib57}
\usepackage{amssymb,amsmath,exscale}
 \numberwithin{equation}{section}

\newcommand{\ZSUno}{\sum _{n=1}^{+\ZIN}}

\newcommand{\ZOMq}{\Omega}

\newcommand{\zg}{\gamma}

\newcommand{\intt}{\int_0^t}
\newcommand{\ints}{\int_0^s}
\newcommand{\intr}{\int_0^r}

\newtheorem{Theorem}{Theorem}

\newcommand{\zaa}{\alpha}

\newcommand{\zt}{\tau}
\newcommand{\zdia}{~~\rule{1mm}{2mm}\par\medskip}
\newcommand{\zthe}{\theta}

\newcommand{\ZLA}{\label}
\newcommand{\ZIN}{\infty}
\newcommand{\zProof}{{\noindent\bf\underbar{Proof}.}\ }

\newcommand{\zzr}{{\rm I\hskip-2.1pt R}}
\newcommand{\ZBI}{\bibitem}
\newcommand{\ZD}{\;\mbox{\rm d}}
\newcommand{\zl}{\lambda}
\newcommand{\ZSI}{\sigma}
\renewcommand{\zthe}{w}

\author{
L. Pandolfi\thanks{Dipartimento di Scienze Matematiche ``Giuseppe Luigi Lagrange'', Politecnico di Torino, Corso Duca degli Abruzzi 24, 10129 Torino, Italy (luciano.pandolfi@polito.it)}
}
\title{Identification of the relaxation kernel in diffusion processes and viscoelasticity with memory via deconvolution\thanks{
This papers fits into the research programme of the GNAMPA-INDAM and has been written in the framework of the   ``Groupement de Recherche en Contr\^ole des EDP entre la France et l'Italie (CONEDP-CNRS)''.}}

\begin{document}

\maketitle
{\bf\underline{Abstract}:}  
 We present an algorithm for the  identification of the relaxation kernel in the theory of diffusion   systems with memory (or of viscoelasticity) which is linear, in the sense that   we propose a linear Volterra integral equation of convolution type whose solution is the relaxation kernel.  The algorithm is based on the observation of the flux through  a part of the boundary of a body. 
 \medskip

\section{Introduction}
The following PDE with persistent memory
\begin{equation}
\ZLA{eq:sistema}
\zthe'=\intt N(t-s)\Delta\zthe(s)\ZD s
\end{equation}
($\Delta $ is the laplacian and the apex denotes time derivative; $w=w(t)=w(x,t)$, whith $x\in\ZOMq$---the variable $x$ is not indicated unless needed for clarity) is encountered in the study of diffusion processes in the presence of complex molecular structures ($\zthe=\zthe(x,t)$ is the concentration, 
see~\cite{PELLE,Hinestroaza}) or in thermodynamics of materials with memory ($\zthe$ is the temperature, 
see~\cite{FabrizioOWENS,PandolfiLIBRO}). It is also encountered in viscoelasticity, usually written   as
\[
\zthe''=N(0)\zthe(x,t)+\intt N'(t-s)\Delta \zthe(s)\ZD s
\] 
(here $\zthe$ is the displacement, see~\cite{PandolfiLIBRO}).

In order to fix the terminology, in this paper we refer to the applications to thermodynamics, and $\zthe$ will be the \emph{temperature.} In this context, we recall that Eq.~(\ref{eq:sistema}) is derived from conservation of energy and a material property for the flux, which takes into account that the flux does not react  instantaneously to the gradient of the temperature, as assumed by the   Fourier law. Rather, it is assumed  that the (density of the) flux at position $x$ and time $t$ is given by
\begin{equation}\ZLA{eq:flux}
q(x,t)=-\intt N(t-s)\nabla \zthe(x,s)\ZD s.
\end{equation}
This material property was first introduced by Maxwell in~\cite{Maxwell} (with  $N(t)=ae^{-bt}$, see also~\cite{cattaneo})  and in the general form by Gurtin and Pipkin in~\cite{GurtinP}.

Dissipation of energy and the second principle of thermodynamics   impose certain restrictions to the \emph{relaxation kernel} $N(t)$
 (see~\cite{Day,GiorgiGentili}) which however is largely unknown and has to be estimated using suitable experiments on samples of the material. For this reason,    identification of $N(t)$ in suitable classes of functions has been widely studied both by engineers and by mathematicians. See   comments and references
  in~\cite{pandKernIdentUNdim}.  We note that $N(t)$ in our model does not depend on $x$ and this is an explicit assumption that the material under examination is homogeneous and isotropic. In this case, $N(t)$ does not depend on the shape of the sample of the material and it can be identified using experiments on samples  with the simplest possible geometry: a bar. Usually, engineering papers (and also our previous paper~\cite{pandKernIdentUNdim}) study this case $
\ZOMq=(a,b)$. It is a fact however that in nondestructive tests, $\zthe=\zthe(x,t)$ where $x\in\ZOMq \subseteq\zzr^d$ and $\ZOMq$ is a bounded region (in the cases of physical interest  $d\leq 3$). In fact, the study of kernel identification in most of the mathematical papers allows $d>1$. 
The methods usually proposed in these papers lead  to the investigation of a system of \emph{ nonlinear} partial differential equations  with memory in the unknowns $\zthe$ and $N$
(see~\cite{FabrizioOWENS,Lerenzilibro}). 

In our previous paper~\cite{pandKernIdentUNdim} we noted that using \emph{two measures} it is possible to identify   the kernel $N(t)$ as the solution of a standard  deconvolution problem, i.e. the solution of a \emph{  linear} Volterra integral equation   of the first kind for the unknown  function $N(t)$. The identification algorithm was proposed in~\cite{pandKernIdentUNdim} in the case $d=1$. Our goal now is the extension of the method to the case $d>1$, a case of interest in nondestructive testing.

It was noted in~\cite{Guidetti}   that, thanks to  the  convolution structure of  system~(\ref{eq:sistema}), if the pair $(N(t), \zthe(t))$  has been identified on an initial  interval $(0,\zt)$ then the identification  of $(N(t),\zthe(t))$  for $t>\zt$ is indeed a linear problem. Our goal can be seen also as  removing the first nonlinear step on $[0,\zt]$ (an interval which in practice cannot be so short to avoid numerical and experimental instabilities). Observe however that the algorithm we propose does not need the explicit computation of $\zthe(t)$.

 Finally we note that the method we are using for kernel identification has been suggested by our study of controllability. A different identification problem whose solution depends on controllability is in~\cite{PandDCDS2}.

\subsection{The assumptions}

The assumptions in this papers concern both $N(t)$ and the region $\ZOMq$. We assume that  $\partial\ZOMq$ is of class $C^2$ (see also the 
comments in Section~\ref{sec:justification}). We single out two (relatively open) subsets $\hat\Gamma$ and $\tilde \Gamma$ of   $\partial\ZOMq$    such that
\[
 {\rm cl}\, \hat \Gamma \cap   
 {\rm cl}\, \tilde \Gamma =\emptyset\,.
\]
 
 The part $\tilde\Gamma$ of the boundary is the part over which the flux  is observed.  
 
As to the kernel $N(t)$ we assume that it is smooth (of class $C^2$) and $N(0)>0$. The restrictions imposed by dissipativity and thermodynamics, in particular $N(t)$ decreasing, can be used for numerical pourposes (see~\cite{pandKernIdentUNdim}) but are not explicitly used in the algorithm. However, we shall see that stability of the null solution  of Eq.~(\ref{eq:sistema}) (which is not an assumption of the identification algorithm) has a practical role in identification too.

 Abel kernels, like $N(t)=at^{-b}$, $a>0$, $b>0$,  are excluded by the regularity assumption.

We contrast these assumptions with   most of the engineering literature, where it is assumed that ${\rm dim}\,\ZOMq=1$ and that the kernel $N(t)$ depends on ``few'' parameters; system~(\ref{eq:sistema}) is numerically solved and then the ``real'' values of the parameters are chosen as the ones which give  a minimum discrepancy from experimental data, see~\cite{BykovMILLISECONDS}.

It is a fact that when $N\in C^2$ then system~(\ref{eq:sistema}) with conditions 
\[
w(x,0)=\xi(x)\,,\quad w(x,t)=f(x,t)\ x\in\partial\ZOMq
\]
is solvable for every $\xi\in L^2(\ZOMq)$ and $f\in L^2(\partial\ZOMq\times (0,T))$ 
and $\zthe\in C([0,T];L^2(\ZOMq))\times C^1([0,T];H^{-1}(\ZOMq))$ (see~\cite{PandolfiLIBRO}).

Finally, let us note a consequence of the assumption $N(0)>0$: variations of the temperature propagates with finite speed $\sqrt{N(0)} $ (in every direction, like in the wave equation, 
see~\cite{DoeschFINITEvelo,PandolfiLIBRO}).

\section{Notations and informal description of the algorithm}

We fix some notations: points of $\partial\ZOMq$ are denoted $x$ and $\ZD\ZSI$ is the surface measure; the points of $\hat\Gamma$ or $\tilde\Gamma$ are denoted respectively $\hat x$ and $\tilde x$ and $\ZD\hat\Gamma$, $\ZD\tilde\Gamma$ the corresponding measures; $\Sigma=\partial\ZOMq\times(0,T)$, $\hat\Sigma=\hat\Gamma\times(0,T)$ $\tilde\Sigma=\tilde\Gamma\times(0,T)$ with measures $\ZD\Sigma$,  $\ZD\hat\Sigma$, $\ZD\tilde\Sigma$.

The exterior normal derivative on $\partial\ZOMq$ is denoted $\zg_1$.

We describe informally the algorithm, which consists in two measures of the flux throughout  $\tilde\Gamma$.  

Every observation is an average against a weighting function $m(x)$ which depends on the instrument.  So, we assume that $m(  x)$ (with support in $\tilde\Gamma$) is given. The observations give    functions
\[
t\mapsto 
\int_{\tilde \Gamma} m(\tilde x) \zg_1 w(\tilde x,t)\ZD\tilde\Gamma\,.
\]

 A first measure is obtained as follows.  We fix $\zthe(x,t)=0$ on $\partial\ZOMq$ and $\zthe(x,0)=\xi(x)$. We measure the corresponding flux on $\tilde\Gamma$.

A second measure is obtained as follows. 
  We fix $\zthe(x,0)$ to a known value (without restriction, $\zthe(x,0)=0$) and we
   impose the temperature $f(\hx,t)$ on $\hat\Gamma$. We choose  
 \begin{equation}\ZLA{eq:formf}
f(\hx,t)=f_0(\hx)f_1(t)\,,\qquad f_1(t)=\intt g(s)\ZD s\,.
 \end{equation}
 Note that the condition $f(\hx,0)=0$ is consistent with $\zthe(x,0)=0$.
  We measure the corresponding flux on $\tilde \Gamma$, using the same instrument, i.e. with the same weight function $m(x)$.
 
  We get respectively
\[
y_\xi( t)= \int_{\tilde\Gamma}  m(\tx) q(\tx,t)\ZD\tilde\Gamma\,,\qquad 
Y^f( t)=\int_{\tilde\Gamma}  m(\tx) q(\tx,t)\ZD\tilde\Gamma \,.
 \]
 It turns out that 
 \[
Y^f( t)=\intt H(t-s)N(s)\ZD s+ \intt g(t-s)y_{\xi_0} (s)\ZD s
 \]
 where $H(t)$  can be computed and $\xi_0$ is a special
 (easily realizable) initial condition so that $y_{\xi_0} $ can be measured. So, the determination of the relaxation kernel $N(t)$ boils down to solve this Volterra integral equation   of the first kind in the unknown $N(t)$, hence to a linear deconvolution problem.
 
 As a last remark, we note that not every instrument for the measure of the flux will do (surely if $m=0$, i.e. the instrument observes nothing and no reconstruction is possible). The  condition on $m(\tilde x)$ is described 
 in~(\ref{RestriSTRUME}).
  
\section{\ZLA{sec:justification}Justification of the algorithm}
 
We introduce the operator $A$ in $L^2(\ZOMq)$, 
\[
{\rm dom} A= H^2(\ZOMq)\cap H^1_0(\ZOMq)\,,\qquad A\phi=\Delta\phi
\]
and an orthonormal basis $\{\phi_n\}$ of $L^2(\ZOMq)$ of eigenvectors, $A\phi_n=-\zl_n^2\phi_n$. 
 We expand the solutions of~(\ref{eq:sistema}) with $\zthe(x,0)=\xi(x)$ and zero boundary temperature is series of the  eigenfunctions $\phi_n$:
\begin{equation}\ZLA{eq:firstWEstudy}
\zthe(x,t)=\ZSUno \phi_n(x)\zthe_n(t)\,,\quad \xi(x)=\ZSUno \xi_n\phi_n(x)\,.
\end{equation}
\paragraph{First we study $y_\xi(t)$.} 
It is easily seen that $\zthe_n(t)$ solves
\[
\zthe_n'(t)=-\zl_n^2\intt N(t-s) \zthe_n(s)\ZD s\,,\qquad \zthe_n(0)=\xi_n\,.
\]
So we have
\[
\zthe_n(t)=  z_n(t)\xi_n
\]
where
\[
z_n'(t)=-\zl_n^2\intt N(t-s) z_n(s)\ZD s\,,\qquad z_n(0)=1 
\]
and
\[
\intt N(t-s)\zthe(x,s)\ZD s=\ZSUno \phi_n(x)\xi_n\intt N(t-s) z_n(s)\ZD s=-\ZSUno \phi_n(x)\frac{1}{\zl_n^2} \xi_n z_n'(t)\,.
\]
The measure of the flux is
\begin{equation} \ZLA{eqFluDatiniz}
  y_{\xi}(t)=\int_{\tilde \Gamma}m(\tx)\left [\ZSUno \zg_1\phi_n(\tx)\frac{1}{\zl_n^2}   z_n'(t)\xi_n\right ]\ZD\tilde\Gamma 
\end{equation}

The series in the bracket  converges in $L^2(0,T; L^2(\tilde\Gamma))$ for every $T>0$ and $\xi\in L^2(\ZOMq)$. The proof is in the Appendix.

\paragraph{Now we study $Y_f(t)$.}We impose the initial condition $\zthe(x,0)=0$ and
we compute the inner product in $L^2(\ZOMq)$ of both the sides of~(\ref{eq:sistema}) with $\phi_n(x)$. 
We use Green formula and we get the following equation for $\zthe_n(t)$ in~(\ref{eq:firstWEstudy}):
 
\[
\zthe_n'(t)=-\zl_n^2\intt N(t-s)\zthe_n(s)\ZD s-\intt N(t-s)\left [\int_{\hat \Gamma}  \zg_1 \phi_n(\hat x) f(\hat x,s)\ZD\hat\Gamma\right ]\ZD s 
\]
and $\zthe_n(0)=0$. We use  $f(x,t)=f_0(\hat x) f_1(t)$ in~(\ref{eq:formf}) and we define

\[
f_{0,n}=\left [\int_{\hat\Gamma}  \zg_1 \phi_n(\hat x) f_0(\hx)\ZD\hat\Gamma\right ]\,.
\] 
Then we have
\[
\zthe_n'(t)=-\zl_n^2\intt N(t-s)\zthe_n(s)\ZD s- f_{0,n} \intt N(t-s) f_1(s)\ZD s
\]
so that (we use $f_1(0)=0$ and the variation of constants formula, see~\cite{PandolfiLIBRO})

\begin{align*}
\zthe_n(t)&= -f_{0,n}\intt z_n(t-r)\intr N(r-s) f_1(s)\ZD s\,\ZD r=\\
&= -f_{0,n}\intt f_1(s)\int_0^{t-s} N(t-s-r) z_n(r)\ZD r\,\ZD s=\\
&=  \frac{1}{\zl_n^2} f_{0,n} \intt  f_1(s)\frac{\ZD}{\ZD s} z_n(t-s)\ZD s=\\
&=  \frac{1}{\zl_n^2} f_{0,n} f_1(t)-\frac{1}{\zl_n^2} f_{0,n} \intt g(s)z_n(t-s)\ZD s \,.
\end{align*}
So we have
\begin{align*}
\zthe(x,t)&=\ZSUno\frac{\phi_n(x)}{\zl_n^2}\left [
f_{0,n} f_1(t)-f_{0,n}\intt g(s) z_n(t-s)\ZD s
\right ]=\\
&=
\left [\left (\ZSUno \phi_n(x) \frac{1}{\zl_n^2}f_{0,n}\right )   f_1(t)-\ZSUno \phi_n(x)\frac{1}{\zl_n^2} f_{0,n} \intt g(s) z_n(t-s)\ZD s\right ]\,.
\end{align*}
It is proved (for example in~\cite{Pandcina,PandolfiLIBRO}) that the series of $w(x,t)$ converges in $C(0,T;L^2(\ZOMq))$.
In order to justify the distribution of the series on the sum, it is sufficient to note that the first series on the right hand side converges in $L^2(\ZOMq)$. This is indeed true since
\[
u_0(x)=\ZSUno\phi_n(x)\frac{1}{\zl_n^2} f_{0,n}
\]
is the solution of 
\begin{equation}\ZLA{probleELLI}
\Delta u=0\ {\rm in}\ \ZOMq\,,\quad
u(x)=f_0^e(x)\  \ x\in \partial\ZOMq\,,\quad f_0^e(x)=\left\{
 \begin{array}{cl} f_0(x)& x\in \hat \Gamma\,, \\ 
0 &x\in  \partial\ZOMq\setminus\hat\Gamma\,.
\end{array}\right.
\end{equation}
 It is known that $ u_0(x)\in H^{1/2}(\ZOMq)\subseteq L^2(\ZOMq)$ and the series converges in this space.

So we have 
\begin{align*}
\nonumber &\intt N(t-s) \zthe(x,s)\ZD s=u(x) \intt N(t-s) f_1(s)\ZD s-\\&-\intt g(s)\left [  \ZSUno \phi_n(x) \frac{1}{\zl_n^2} f_{0,n} \int_0^{t-s} N(t-s-r) z_n(r)\ZD r  \right ]\ZD s=\\
 &=u(x) \intt N(t-s) f_1(s)\ZD s+\intt g(t-s) \left [
\ZSUno \phi_n(x)\frac{1}{\zl_n^4} f_{0,n}z_n'(s)
\right ]\ZD s\,.
\end{align*}

So we have:
\begin{align}
\nonumber Y^f(t)&=\left [\int _{\tilde\Gamma } m(\tilde x)  \zg_1 u(\tilde x)\ZD\tilde\Gamma\right ]\intt N(t-s)f_1(s)\ZD s+\\
\ZLA{eq:preliflu}&+\intt g(t-s)\int _{\tilde\Gamma}m(\tilde x)\left [\ZSUno \zg_1\phi_n(\tilde x)\frac{1}{\zl_n^4}f_{0,n}z_n'(s)\right ]\ZD\tilde\Gamma\, \ZD s\,.
\end{align}
The convergence in $L^2(0,T;L^2(\tilde\Gamma))$ of the series in  the bracket is in the appendix.

We compare the last bracket with~(\ref{eqFluDatiniz}) and we see that
\[
\intt g(t-s)\int _{\tilde\Gamma}m(\tilde x)\left [\ZSUno \zg_1\phi_n(\tilde x)\frac{1}{\zl_n^4}f_{0,n}z_n'(s)\right ]\ZD\tilde\Gamma\, \ZD s=\intt g(t-s) y_{\xi_0}(s)\ZD s 
\]
where
\[
 \xi_0(x)=\ZSUno\phi_n(x) \frac{1}{\zl_n^2}f_{0,n}=u_0(x)\,.
\]
So, once the measure $y_{\xi_0}(s)$ is known, identification of $N(t)$ is reduced to a standard deconvolution problem,
provided that we can give a meaning to
\begin{equation}\ZLA{eq:osseU}
\left [\int _{\tilde\Gamma } m(\tilde x)  \zg_1 u_0(\tilde x)\ZD\tilde\Gamma\right ] 
\end{equation}
and provided that
\begin{equation}\ZLA{RestriSTRUME}
\left [\int _{\tilde\Gamma } m(\tilde x)  \zg_1 u_0(\tilde x)\ZD\tilde\Gamma\right ]\neq 0\,.
\end{equation}
This last propery is a property of the instrument used for the  measure and we assume that it is satisfied.

\subparagraph{We comment the conditions under which the integral~(\ref{eq:osseU}) makes sense.} This depends on the smoothness of $\partial\ZOMq$,  of the weight $m(\tilde x)$ and of $f_0^ e(x)$.

When $f^e _0\in L^2(\ZOMq)$ ( defined in~(\ref{probleELLI})   then we have $u_0\in H^{1/2}(\ZOMq)\subseteq L^2(\ZOMq)$ and $\Delta u_0=0\in L^2(\ZOMq)$. Then we have $\zg_1 u_0\in H^{-3/2}(\Gamma)$ (see~\cite[Ch.~2]{LionsMagenes}) and the  bracket is well defined if $m(x)$ is   smooth, $m\in H^{3/2}(\ZOMq)$. This also in the case that $\hat\Gamma$ and $\tilde\Gamma$ do intersect. But, we can also take advantage of the fact that  ${\rm cl} \hat\Gamma \cap{\rm cl}  \tilde\Gamma =\emptyset$.

We assume more regularity of $\partial\ZOMq$, i.e. $\partial\ZOMq\in C^{2,\zaa}$ (any $\zaa>0$). If $f^e_0$ is bounded, problem~(\ref{probleELLI}) admits a  solution  which is continuous on $\tilde\Gamma$ (see~\cite[Lemma~2.13]{Gilbargtrudinger}) (note that $f^e_0$ is zero, hence continuous on $\tilde\Gamma$ because ${\rm cl}\, \hat\Gamma \cap {\rm cl}\,\tilde\Gamma=\emptyset$).
Moreover,  it follows from~\cite[Lemma~6.18]{Gilbargtrudinger} that $\zg_1 u(\tilde x)$ is continuous on $\tilde\Gamma$ and~(\ref{eq:osseU}) makes sense for every $m_0\in L^2(\tilde\Gamma)$, in  particular if it is the characteristic function of a subset of $\tilde\Gamma$. 

A remark is in order here:~\cite[Lemma~6.18]{Gilbargtrudinger} assumes $u\in C({\rm cl}\,\ZOMq)$ but, as simply seen from the proof and explicitly stated at~\cite[p.~106]{Gilbargtrudinger}, the result holds also if $u$ is continuous on $\ZOMq_1$ where $\ZOMq_1\subseteq\ZOMq$ is such that $\tilde\Gamma\subseteq\partial \ZOMq_1$.  The region  $\ZOMq_1$ exists  because $ {\rm cl}\,\hat \Gamma \cap {\rm cl}\,\tilde \Gamma =\emptyset $.

\subparagraph{On the initial condition $\xi_0=u_0$.}
A final point has to be settled. While it is technically reasonable to impose   boundary conditions which are more or less arbitrary,   it is not such an easy task to impose an initial condition unless it has some special property. In fact, the material has to be prepared to a certain temperature distribution and then released  to evolve freely for $t $ larger then an ``initial'' time, without restriction for $t>0$. This is in general difficult.
But, we are interested here to the \emph{special} initial condition $\xi_0(x)=u_0(x)$ and this initial condition can be easily achieved \emph{if system~(\ref{eq:sistema}) is stable,}  i.e. if
 for every $\xi$ the corresponding solution (with zero boundary temperature) $\zthe(t)\to 0$ in $L^2(\ZOMq)$, for $t\to +\ZIN$ (as it usually happens in practice since   the memory kernel $N(t)$ dissipates enough energy).
In this case the initial condition $\xi_0$ can easily be imposed since $\xi_0=u_0$ solves~(\ref{probleELLI}). To see this,
let us consider Eq.~(\ref{eq:sistema}) with $\zthe(x,0)=0$ and
\[
\zthe(x,t)=f_0^e(x)\quad x\in\partial\ZOMq
\]
(note that the boundary condition is constant in time). 

The function
\[
\eta(x,t)=\zthe(x,t)-\xi_0(x)=\zthe(x,t)-u_0(x)
\]
solves Eq.~(\ref{eq:sistema}) with the boundary condition put equal zero and the initial condition equal to $\xi_0$.

Stability implies $\lim_{t \to+\ZIN} \eta(\cdot,t)=0 $ in $L^2(\ZOMq)$ so that
\[
\lim_{t \to+\ZIN}\zthe(x,t)=\xi_0(x)\,:
\]
in order to realize the special initial condition $\xi_0$ it is sufficient to apply the boundary temperature $f_0^e$ (constant in time) to the body, for a time large enough, after which the heated external source is removed and  the measure $y_{\xi_0}(t)$ can be taken.
 
\section{Appendix:  the series~(\ref{eqFluDatiniz})}
 
We need few information. A sequence $\{e_n\}$ in a (real or complex) Hilbert space is a \emph{Riesz sequence} if it can be transformed to an orthonormal basis (possibly of a   Hilbert space of different dimension) using a linear bounded and boundedly invertible transformation (see~\cite{Heil,PandolfiLIBRO} for Riesz sequences). It is a fact that when $\{e_n\}$ is a Riesz sequence then the series 
 
\begin{equation}\ZLA{eq:seriePERconv}
\ZSUno h_n e_n
\end{equation}
converges  if and only if $\{h_n\}\in l^2$ and the sum does not depend on the order of the elements of the series. 
 
Let $\{e_n\}$ be a sequence in a Hilbert space $H$. We associate to this sequence the \emph{moment} operator $\Mo$: $H\mapsto l^2$ 
\[
\Mo h=  \{ \langle h,e_n\rangle\} 
\] 
(here   $\langle\cdot,\cdot\rangle$ denotes the inner product in $H$).
 \emph{The sequence $\{e_n\}$ is a Riesz sequence if and only if the operator $\Mo$ is continuous and surjective.}

We recall our goal, which is the proof of the convergence   in $L^2(\tilde\Gamma\times(0,T))$ of the series in~(\ref{eqFluDatiniz}), i.e.
\begin{equation}\ZLA{eq:Da provaCONVE}
\ZSUno \xi_n\left (\zg_1\phi_n\right ) \intt N(t-s) z_n(s)\ZD s\,.
\end{equation}
We prove convergence in $L^2(\partial\ZOMq\times (0,T))$ for every $T>0$.
This property depends on  the fact that system~(\ref{eq:sistema}) is controllable in $L^2(\ZOMq)$ using controls $f\in L^2(\partial\ZOMq\times(0,T))$ for $T$ large enough.
Namely, there exists a time $T_0$ such that for every $\xi\in L^2(\ZOMq)$ there exists $f\in  L^2(\partial\ZOMq\times(0,T_0))$ such that $\zthe(T_0)=\xi$ (the initial condition is zero and we assume $\zthe=f$ on $\partial\ZOMq$). Using similar arguments as in the study of $Y_f(t)$, we expand $\zthe(x,t)$ in series of the eigenfunctions $\phi_n(x)$,
\[
\zthe(x,t)=\ZSUno \zthe_n(t)\phi_n(x)\,.
\]
 The functions $\zthe_n(t)$ are given by
\[
\zthe'_n(t)=-\zl_n^2\intt N(t-s)\zthe_n(s)\ZD s-\int_0^t N(t-s) \int_{\partial\ZOMq} f(x,s)\zg_1\phi_n\ZD\ZSI\,\ZD s\,.
\]
Hence
 
\[
\zthe_n(t)=-\intt z_n(t-\zt)\left [ \int_0^\zt N(\zt-s) \int_{\partial\ZOMq} f(x,s)\zg_1\phi_n\ZD\ZSI\,\ZD s\right ]\ZD \zt\,.
\]
Analogously, we expand $\xi=\ZSUno\xi_n\phi_n(x)$ and we see that controllability at time $T_0$ is equivalent to the solvability  of the following  {moment problem}:
\begin{align*}
\Mo f &= \int_0^{T_0} z_n(T_0-\zt)\left [ \int_0^\zt N(\zt-s) \int_{\partial\ZOMq} f(x,s)\zg_1\phi_n\ZD\ZSI\,\ZD s\right ]\ZD \zt=\\
&= \int_0^{T_0} \int_{\partial\ZOMq} f(x,T_0-t)\left [
\zg_1\phi_n\ints N(s-r) z_n(r)\ZD r\ZD\ZSI
\right ] =
\xi_n\,.
\end{align*}
 Here $\{\xi_n\}$ is an arbitrary $l^2$ sequence.
 
The fact that $f\mapsto \zthe $ is a linear continuous transformation from $L^2(\partial\ZOMq\times(0,T_0))$ to $C([0,T_0];L^2(\ZOMq))$  implies that the {moment operator} $\Mo$ is continuous from 
$L^2(\partial\ZOMq\times(0,T_0))$ to $l^2$, and controllability shows that it is surjective. So, the sequence of functions
\[
\left \{  
\zg_1\phi_n\ints N(s-r) z_n(r)\ZD r\ZD\ZSI
 \right \}
\]
is a Riesz sequence in $L^2(\partial\ZOMq\times (0,T_0))$ and this in turn implies that the series 
\[
\ZSUno \xi_n  \left [
\zg_1\phi_n\ints N(s-r) z_n(r)\ZD r\ZD\ZSI
\right ]
\]
 converges in $L^2(\partial\ZOMq\times (0,T))$ for every $T$ (both larger and and less then $T_0$) and for every $\{\xi_n\}\in l^2$.
 
In particular, it converges in $L^2(\tilde\Gamma\times (0,T))$ and this argument applies in particular to the series~(\ref{eq:preliflu}). 

Note that once convergence has been proved, formula~(\ref{eqFluDatiniz}) can be written as
\begin{align*}
 &\ZSUno \frac{\xi_n}{\zl_n^2}\left (
 \int _{\tilde\Gamma} m(\tx)\left (\zg_1\phi_n(\tx)\right )\ZD\tilde\Gamma
 \right ) z_n'(t)=\\
&=
 \ZSUno \xi_n \left (
  \int _{\tilde\Gamma} m(\tx) \zg_1\phi_n(\tx)   \ZD\tilde\Gamma
 \right ) \zl_n\intt N(t-s) z_n(s)\ZD s\,.
 \end{align*}
 
 Expecially for numerical reasons, it makes sense to elaborate on formula~({\ref{eq:Da provaCONVE}). It is known that each eigenvalue of $A$ has finite multiplicity. If we denote $\{-\zl_n^2\}$ the sequence of the \emph{distinct} eigenvalues of $A$,  $\{\phi_{n,k}\} _{k=1}^{k_n} $  will be the eigenvectors with the same eigenvalue $-\zl^2_n$.
 We note that $z_n$ depends only on $-\zl_n^2$ and not on the eigenfunction, so that~(\ref {eq:Da provaCONVE}) can be written as
 \[
\ZSUno \left (\sum _{k=1}^{k_n} \xi _{n,k} \Psi_{nk}\right ) \left [\zl_n\intt N(t-s) z_n(s)\ZD s\right ]\,,\qquad \Psi_{n,k}=\frac{1}{\zl_n}\zg_1\phi _{n,k}\,. 
 \]
 We prove:
 \begin{Theorem}
 The sequence $\left\{  \left \|    \sum _{k=1}^{k_n} \xi _{n,k}\Psi_{n,k}   \right \|_{L^2(\partial\ZOMq)} \right \} $ belongs to $l^2$.
 \end{Theorem}
 \zProof If $k_n=1$ for every $n$ then the result follows since $\{\Psi_n\}=\{\Psi_{n,1}\}$ is bounded in $L^2(\partial\ZOMq)$, see~\cite{TAOcorrect}. 
 We consider the wave equation 
 
 \begin{equation}\ZLA{eq:Wave}
u''=\Delta u\,,\quad u(x,0)=0\,, \ u'(x,0)= \xi\,,\qquad u=0\ {\rm on}\ \partial\ZOMq\,.
 \end{equation} 
 It is known that the solution $u$ of this equation satisfies
\begin{equation}
 \ZLA{eq:directINEQ}
 \| \zg_1 u\|^2_{L^2(\Sigma)} \leq M\|\xi\|^2_{L^2(\ZOMq)}
 \end{equation}
 (the constant $M$ does not depend on $\xi$, see~\cite{LIONSlibro}).
 
 We apply   inequality~(\ref{eq:directINEQ}) to the solution of~(\ref{eq:Wave}) when
 \[
\xi= \sum _{k=1}^{k_n} \xi_{n,k}\phi_{n,k}\,. 
 \]
 The solution is
 \[
u(x,t)=\left (\frac{\sin \zl_n t}{\zl_n}  \right ) \sum _{k=1}^{k_n} \xi_{n,k}\phi_{n,k}
 \]
and inequality~(\ref{eq:directINEQ}) gives
\begin{align*}
\left (\int _{0}^T \sin^2\zl_n t\ZD t\right )
 \int _{\partial\ZOMq}\left \|\sum _{k=1}^{k_n} \xi_{n,k}\Psi_{n,k}\right \|^2\ZD\sigma\leq 
   M\sum _{k=1}^{k_n} |\xi_{n,k}|^2
    \ {\rm   i.e.}\\
  \left \|\sum _{k=1}^{k_n} \xi_{n,k}\Psi_{n,k}\right
  \|^2_{L^2(\partial\ZOMq)}
    \leq M \sum _{k=1}^{k_n} |\xi_{n,k}  | ^2
\end{align*}
(for a different constant $M$). 

The constant does not depend on $n$ and so, summing on $n$, we find
\[
\ZSUno \left \|\sum _{k=1}^{k_n} \xi_{n,k}\Psi_{n,k}\right \|^2_{L^2(\partial\ZOMq)}\leq M\|\xi\|^2_{L^2(\ZOMq)}\,.
\]
This is the required inequality.\zdia
  

\begin{thebibliography}{99}

  \ZBI{FabrizioOWENS} Amendola, G. Fabrizio, M. Golden, J.M.,  Thermodynamics of materials with memory. Theory and applications. Springer, New York, 2012. 

\ZBI{PELLE} Barbeiro, S., Ferreira J.A., Coupled vehicol-skin models for drug release, \emph{Comp. Methods Appl. Mech. Engrg.} {\bf 198} 2078-2086 (2009).


\ZBI{BykovMILLISECONDS}       Bykov. D.L., Kazakov, A.V., Konovalov, D.N., Me'lnikov, V.P., Osavchuk, A.N., Peleshko, V.A.,     Identification of the model of nonlinear viscoelasticity of filled polymer materials in millisecond time range, Mechanics of Solids, 47 641-645, 2012     


\ZBI{cattaneo}  Cattaneo, C.:  {Sulla conduzione del calore.} \emph{Atti   Semin. Mat. e Fis.  Univ. di Modena}~\textbf{3},
   83-101 (1948)
   
\ZBI{Day} Day, W.A., On monotonicity of the relaxation functions of viscoelastic materials, \emph{Proc. Camb. Phil. Soc.} {\bf 67} 503-508 1970   
   
\ZBI{Hinestroaza}  De Kee, D.,     Liu, Q.,   Hinestroza, J.: Viscoelastic (non-fickian) diffusion. {\em The Canada J. of chemical engineering\/}~\textbf{83},   913-929 (2005).

\ZBI{DoeschFINITEvelo}  Desch, W.,  Grimmer, R.C.:  Initial-boundary value problems for integro-differential equations.    \emph{J. Integral Equations}~\textbf{10},   73-97 (1985) 

\ZBI{GiorgiGentili} Giorgi, C., Gentili, G., Thermodynamics properties and stability for the heat flux equation with linear memory, \emph{Quarterly Appl. Math.,} {\bf 51}  343-362, 1993.

  \ZBI{Gilbargtrudinger} Gilbarg, D. Trudinger, N.S., Elliptic partial differential equations of second order, Springer-Verlag, Berlin, 1977.
  \ZBI{Heil} Heil, C.A., Basis theory primer. Birkh\"auser Springer, New York, 2011.
  

  \ZBI{Guidetti} Guidetti, D.,  Reconstruction of  variation convolution kernel in an abstract wave equation. \emph{Forum Math.} {\bf 6} 1129-1160, 2010. 

\ZBI{GurtinP} Gurtin, M.E.,   Pipkin, A.G.:  
  A general theory of heat conduction with
 finite wave speed. \emph{
 Arch. Rational Mech. Anal.}~\textbf{31}, 113-126 (1968)
 
 \ZBI{TAOcorrect}  A. Hassel, T. Tao, Erratum for ``Upper and lower bounds for normal derivatives of Dirichlet eigenfunctions''. {\em Math. Res. Lett.\/} {\bf 17} (2010)  793--794.
 
 
 
 \ZBI{LIONSlibro} Lions, J-L.,   \emph{Contr\^olabilit\`e exacte, perturbations et stabilization de syst\'emes distribu\`es.\/}
Vol.~1, Recherches en Math\'e matiques Appliqu\'e e~9, Masson, Paris (1988)

\ZBI{LionsMagenes} Lions, J.-L., Magenes, E., \emph{Problèmes aux limites non homogènes et applications,} Vol. 1.  Dunod, Paris 1968

\bibitem{Lerenzilibro} Lorenzi, A., An introduction to identification problems via functional analysis, VSP, Utrecht,  2001.

 

\ZBI{Maxwell} Maxwell, J.C.,  On the dynamical theory of gases, \emph{Philos. Trans.,} \textbf{157} 49-88   1867

 
 
  
   \bibitem{Pandcina}  Pandolfi, L.: Riesz systems and moment method in the study of heat equations with memory  in one space dimension.  \emph{Discrete Contin. Dyn. Syst. Ser. B.\/}~\textbf{14}, 1487-1510 (2010)
   
 
   
     \bibitem{PandDCDS2}  
 Pandolfi, L.:
 { Riesz systems and an identification problem for heat equations with memory.\/} 
 \emph{Discrete Contin. Dyn. Syst. Ser. S}~\textbf{4},  745--759 (2011)
 
 
  
 
 \ZBI{PandolfiLIBRO}  Pandolfi, L.,
Distributed systems with persistent memory. 
Control and moment problems. Springer Briefs in Electrical and Computer Engineering. Control, Automation and Robotics. Springer, Cham, 2014.

\ZBI{pandKernIdentUNdim} Pandolfi, L., A linear algorithm for the identification of a relaxation kernel using  two boundary measures \emph{Inverse Problems,} {\bf 31} no. 10, 105003, 2015. 
 
  \end{thebibliography}
 \enddocument